\documentclass{article}
\usepackage{psfig}

\def\qed{\hfill {\hbox{${\vcenter{\vbox{               
   \hrule height 0.4pt\hbox{\vrule width 0.4pt height 6pt
   \kern5pt\vrule width 0.4pt}\hrule height 0.4pt}}}$}}}

\newtheorem{theorem}{Theorem}
\newtheorem{definition}{Definition}

{\qed\par\medskip}

\begin{document}

\markboth{Louis H. Kauffman and Slavik V. Jablan} {A theorem on
amphicheiral alternating links}

\centerline{\bf A NOTE ON AMPHICHEIRAL
ALTERNATING LINKS}

\bigskip 

\centerline{\footnotesize LOUIS H. KAUFFMAN, SLAVIK V. JABLAN$^ \dag
$}

\centerline{\footnotesize\it University of Illinois at Chicago,
Department of Mathematics,}\centerline{\footnotesize\it Statistics
and Computer Science (m/c 249),}\centerline{\footnotesize\it 851
South Morgan Street, Chicago,}\centerline{\footnotesize\it Illinois
60607-7045, USA} \centerline{\footnotesize\it kauffman@uic.edu}

\centerline{\footnotesize\it The Mathematical Institute$^\dag$, Knez
Mihailova 36,}\centerline{\footnotesize\it P.O.Box 367, 11001
Belgrade, Serbia} \centerline{\footnotesize\it sjablan@gmail.com}

\medskip



\bigskip

\begin{abstract}
\noindent After reconsidering the Dasbach-Hougardy counterexample to the
Kauffman Conjecture on alternating knots, we reformulate the
conjecture and consider Dasbach-Hougardy counterexample and similar counterexamples
in the light of the reformulated conjecture.
\end{abstract}

\section{Introduction}
\begin{figure}[th]
\centerline{\psfig{file=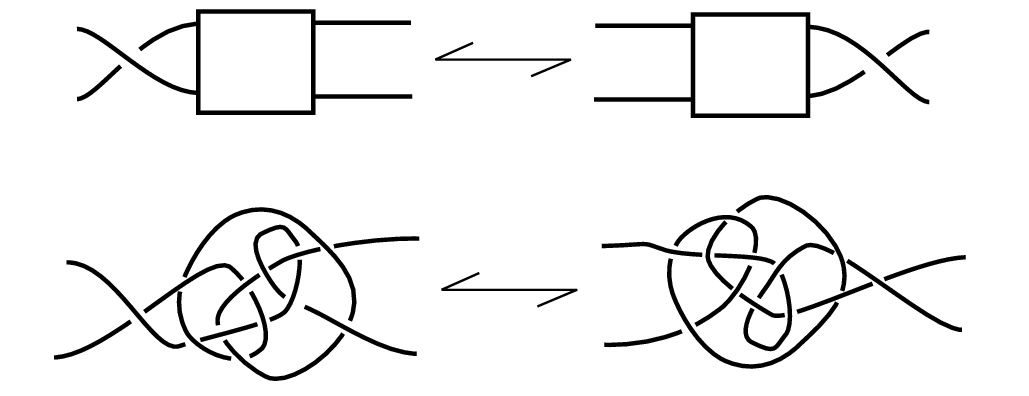,width=7cm}} \vspace*{8pt}
\caption{Flyping.  \label{flype}}
\end{figure}

We begin with a short description of the general background to the paper.
Whitney \cite{9} characterized spherical projections of planar graphs. He showed that two embeddings
of a planar graph into a sphere are related by $2$-isomorphisms. A result similar in sprit in knot theory
is the Flyping Theorem that was conjectured by Tait and proven by Menasco and Thistlethwaite: any
two reduced alternating projections of a prime alternating link are related by {\it flypes}. A flype is
a replacement of a subtangle in the knot diagram as shown in Figure~\ref{flype}. The replacement
corresponds to turning the tangle around by $180$ degrees and pulling its connecting arcs as well
so that a crossing on one side of the tangle ends up on the other side, as shown in Figure 1.
Flyping is an ambient isotopy, and if a link is alternating and one performs a flype, then the resulting link is also alternating.

\begin{figure}[th]
\centerline{\psfig{file=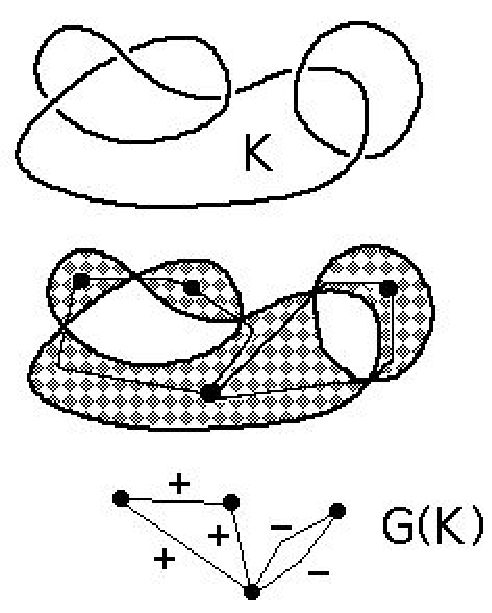,width=6cm}} \vspace*{8pt}
\caption{The Checkerboard Graph $G(K).$ \label{checker}}
\end{figure}

In Figure~\ref{checker} we illustrate the formation of the {\it checkerboard graph} $G(K)$ of a link diagram $K.$ The checkerboard graph of a link diagram $K$ in the plane is formed by first coloring the
regions of the link diagram (shaded and unshaded) with two colors so that adjacent regions have
distinct colors, and so that the outer unbounded region is unshaded. Then a node is assigned to each
shaded region, and edges occur in the graph whenever two regions share a crossing in the link
diagram. In the Figure~\ref{checker} we
have also illustrated how to decorate the graph $G(K)$ with signs so that the diagram $K$ can be
reconstructed from it. In this paper we shall not need to keep track of these signs and so we denote
by $G(K)$ the undecorated graph that is associated with the checkerboard coloring. The reader will find it easy to verify that the planar dual $G^{*}(K)$ of the graph $G(K)$ is obtained by forming a graph
by the same prescription but  using the unshaded regions.  We let $K^{*}$ denote the mirror image of
the diagram $K$ that is obtained by switching all the crossings of $K.$ It is not hard to see that
$G(K^{*})$ is graph isomorphic {\it on the surface of the two dimensional sphere}  with the dual graph $G^{*}(K).$ This leads to natural questions
about the relationship of graphs and dual graphs for links  that are {\it amphichieral.} A link is said
to be amphichieral if it is ambient isotopic to its mirror image. It should be remarked that the natural domain for graph isomorphisms in this paper is on the surface of a two dimensional sphere. Knots and
links are represented by planar diagram embeddings so that the checkerboard graphs are well-defined.

For denoting knots and links we use Conway notation \cite{1,4}.
All tangles are finite compositions of elementary tangles. The
elementary tangles are 0, 1 and $-1$. Algebraic knots and links can be obtained by using three operations: {\it sum}, {\it product}, and {\it ramification},
leading from tangles $a$, $b$ to the new tangles $a+b$, $-a$, $a\,b$, and $(a,b)$, respectively, where $-a$ is the
image of $a$ in NW-SE mirror line, $a\,b = -a+b$, $-a = a\,0$, and $(a,b) = -a-b$ (Fig. 3). Polyhedral knots and links can be obtained by substituting vertices of basic polyhedra by algebraic tangles.

\begin{figure}[th]
\centerline{\psfig{file=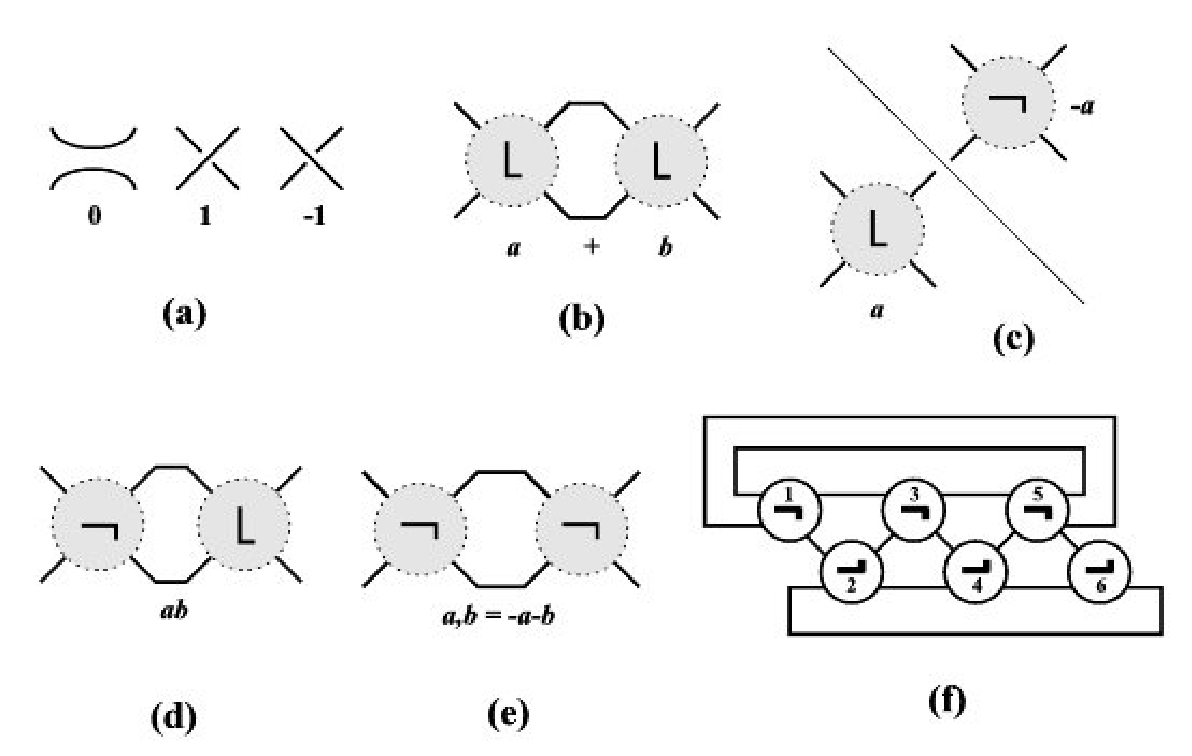,width=4.20in}} \vspace*{8pt}
\caption{(a) The elementary tangles; (b) sum of tangles $a+b$; (c) tangle $-a = a\,0$; (d) product of tangles $a\,b = -a+b$; (e) ramification of tangles $(a,b) = -a-b$; (f) the basic polyhedron $6^*.$  \label{tangles}}
\end{figure}

Louis Kauffman conjectured in \cite{6} (revised in \cite{7}) that every
amphicheiral alternating knot can be drawn so that the checkerboard
graph of the knot diagram is self dual:

\noindent {\bf Conjecture 1.1.} {\it (Kauffman Conjecture) Let $K$
be an alternating amphicheiral knot. Then there exists a reduced
alternating knot diagram $D$ of $K$, such that $G(D)$ is isomorphic
to $G^*(D)$, where $G(D)$ is a checkerboard-graph of $D$ and
$G^*(D)$ is its dual \cite{7}.}

In the statement of the Kauffman Conjecture the term "isomorphic" refers to isomorphism of abstract graphs, denoted in this paper by $\simeq $.

In the paper \cite{2} Oliver Dasbach and Stephan Hougardy found a
counterexample to this conjecture,  an amphicheiral alternating knot $14a_{10435}$, given
in Conway notation \cite{1,4} as $(2\,1,3)\,1\,1\,(2\,1,3)$, which has four
different minimal diagrams $D_1$-$D_4$ (Fig.
4 a,b,c,d). None of their corresponding checker-board graphs $G_1$, $G_2$,
$G_3$, $G_4$ is isomorphic to its dual $G^*_1$, $G^*_2$, $G^*_3$,
$G^*_4$. However, $G_3\simeq G^*_1$ (Fig. 5a) and $G_4\simeq G^*_2$
(Fig. 5b).

\begin{figure}[th]
\centerline{\psfig{file=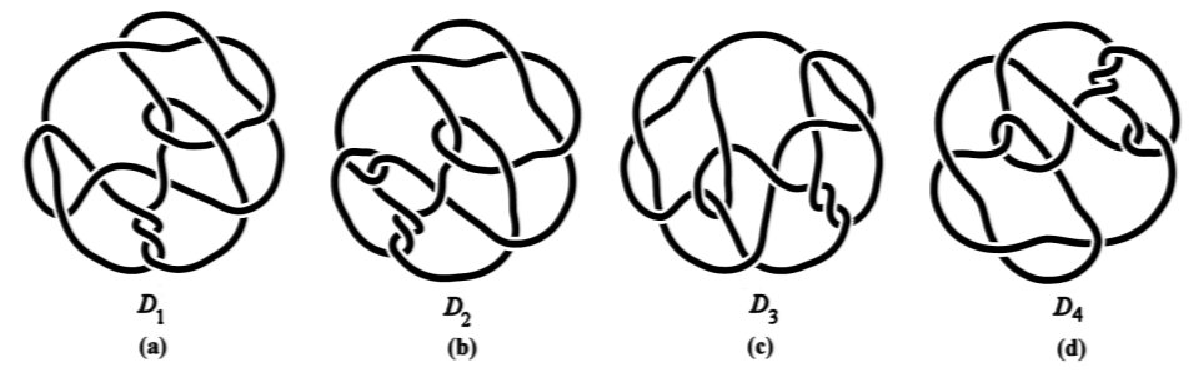,width=4.50in}} \vspace*{8pt}
\caption{Four minimal diagrams $D_1$-$D_4$ of the knot
$(2\,1,3)\,1\,1\,(2\,1,3)$. \label{f1.1}}
\end{figure}

\begin{figure}[th]
\centerline{\psfig{file=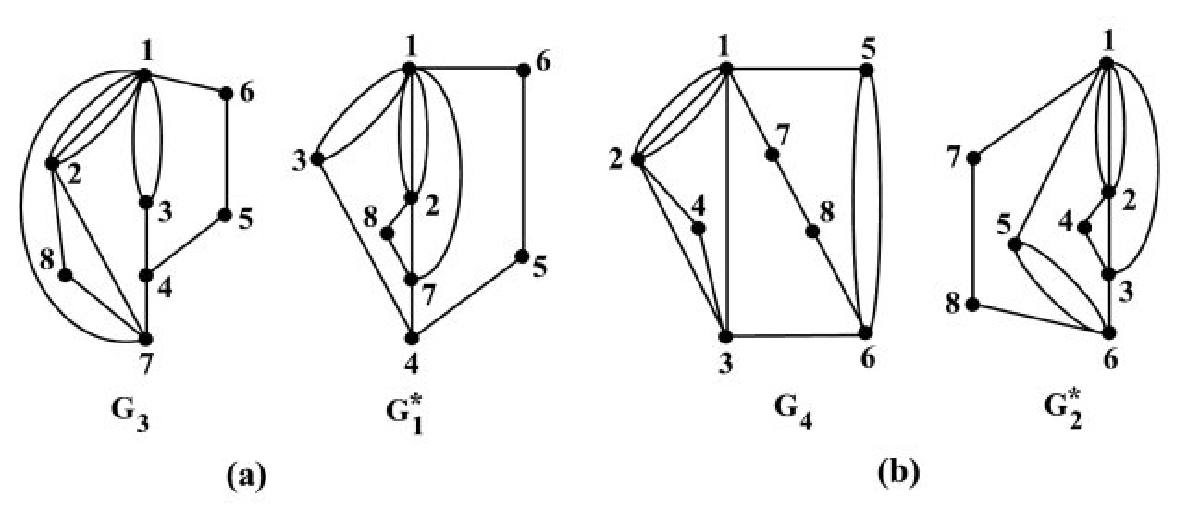,width=4.50in}} \vspace*{8pt}
\caption{Isomorphic graphs (a) $G_3\simeq G^*_1$; (b) $G_4\simeq
G^*_2$ . \label{f1.2}}
\end{figure}

\section{Reformulated Kauffman Conjecture}
This motivates us to reformulate the original Kauffman conjecture for alternating
knots  using the Tait Flyping Conjecture, proved by
Menasaco and Thistlethwaite \cite{8}. This result tells us that two reduced alternating
link diagrams are ambient isotopic if and only if they are related by a series of flypes
as illustrated in Figure 1. Given a knot or link $K$, let $G(K)$ denote the checkerboard
graph for $K$. We will say that that two graphs $G(K)$ and $G(K')$ are flype equivalent
if $K$ and $K'$ can be transformed into one another by a sequence of flypes.

Now suppose that $K$ is a reduced alternating diagram that is ambient isotopic
to its mirror image diagram $K^*$. Let $G(K)$ denote the graph of $K$ and $G(K^*)$
denote the graph of the mirror image knot $K^*$. Then we know from construction
that $G(K^*) = G^*(K)$ where $G^*$ denotes the dual graph of $G(K)$ in the plane
or on the surface of the two-dimensional sphere. Since we know that $K^*$ is related
to $K$ by flypes (by the theorem of Menasaco and Thistlethwaite),
it follows that the graphs $G(K)$ and $G^*$ are flype equivalent.
This is the correct statement that can replace the original Kauffman conjecture.

In other words, we have the following correct statement that replaces the original Kauffman Conjecture:
\bigbreak

\begin{theorem}
Let $K$ be an reduced alternating, prime,
knot or link. Then the checkerboard graph $G(K)$ is flype equivalent to
its dual graph $G^{*}(K)$ if and only if $K$ is amphicheiral.
\end{theorem}

\noindent Note that this statement implies that if $K$ is amphicheiral, then some minimal diagram of $K^{*}$  has checkerboard graph that is equivalent to
the dual of the checkerboard graph of $K.$

The interest in this simple reformulation, in the
light of the Flyping Theorem, is that if we are given a minimal diagram for an
alternating link, then this diagram is known \cite{5} to be alternating and it is not hard
to enumerate all the possible flype-equivalent diagrams.
\bigbreak

\noindent {\bf Remark.} Note that one can take the
chiral 3-component link with $n=16$ crossings
$6^*(2\,1,2)\,1.(2,2\,1)\,1$, with 16 different minimal diagrams.
Its minimal diagrams $D_1$  (Fig. 6a) and $D_2$  (Fig. 6b) satisfy the
relationship $G_1\simeq G^*_2$ (Fig. 7a,b). This shows that some dual graphs may be
equivalent to some graphs of the original link but that this does not, in itself, imply achirality or amphicheirality.

\begin{figure}[th]
\centerline{\psfig{file=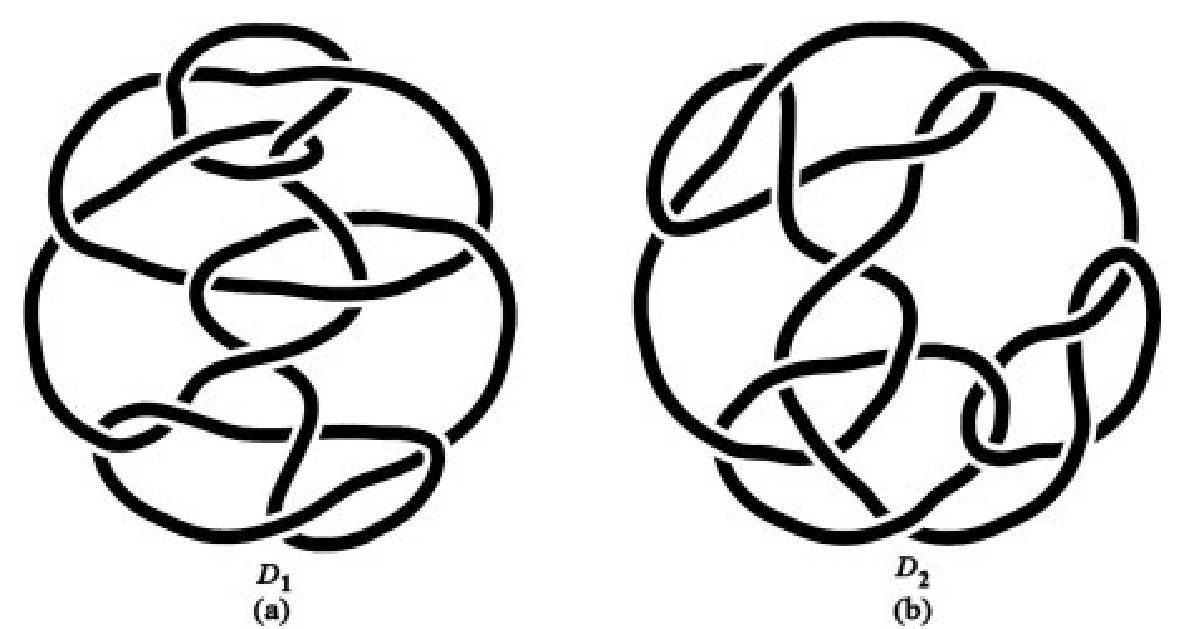,width=2.50in}} \vspace*{8pt}
\caption{Two minimal diagrams (a) $D_1$;
(b) $D_2$ of the chiral link
$6^*(2\,1,2)\,1.(2,2\,1)\,1$ with $G_1\simeq G^*_2$. \label{f1.3}}
\end{figure}

\begin{figure}[th]
\centerline{\psfig{file=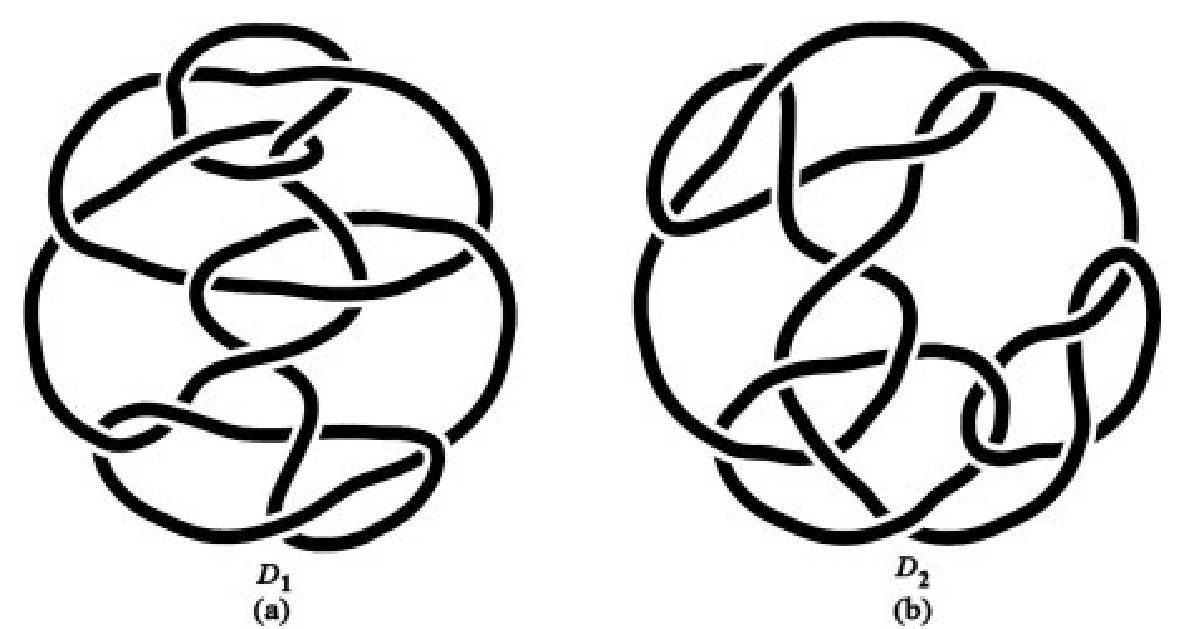,width=2.50in}} \vspace*{8pt}
\caption{Two minimal diagrams (a) $D_1$;
(b) $D_2$ of the chiral link
$6^*(2\,1,2)\,1.(2,2\,1)\,1$ with $G_1\simeq G^*_2$. \label{f1.3}}
\end{figure}

\begin{figure}[th]
\centerline{\psfig{file=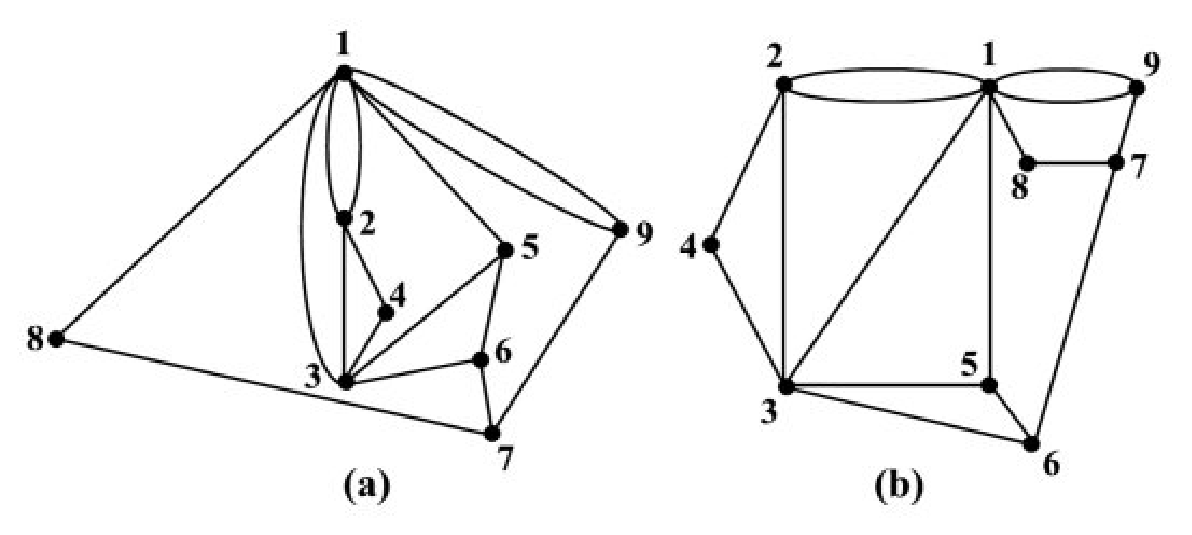,width=2.20in}} \vspace*{8pt}
\caption{Isomorphic graphs $G_1$ and $G^*_2$.
\label{f1.4}}
\end{figure}

For links with a single minimal diagram, the reformulated Kauffman Conjecture can be
used as the criterion for recognition of amphicheiral alternating
links: a link with a single minimal diagram $D$ is amphicheiral {\it
iff} $G(D)\simeq G^*(D)$. For example, among mutant alternating
knots $K_1=.(2,3).(3,2)$ and $K_2=.(2,3).(2,3)$ with a single
minimal diagram (Fig. 8) the first is amphicheiral, and the other is
not, because the graph of the first knot is self-dual, and the other
is not (Fig. 9).

\begin{figure}[th]
\centerline{\psfig{file=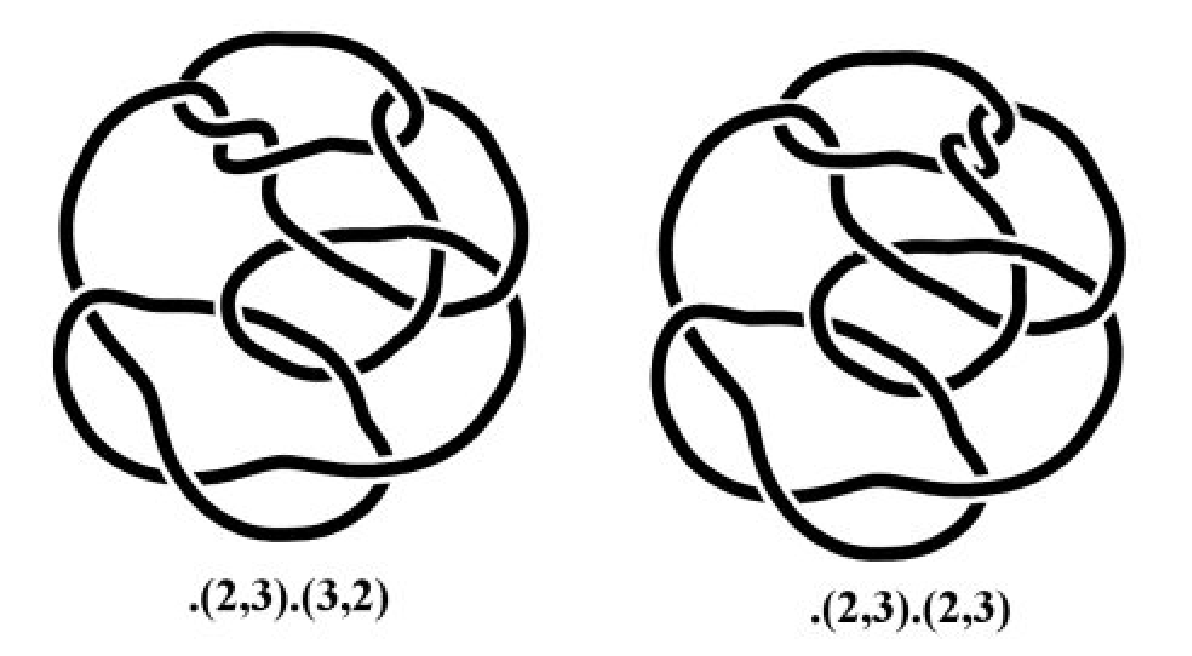,width=2.50in}}
\vspace*{8pt} \caption{Knots $K_1=.(2,3).(3,2)$ and
$K_2=.(2,3).(2,3)$. \label{f1.5}}
\end{figure}

\begin{figure}[th]
\centerline{\psfig{file=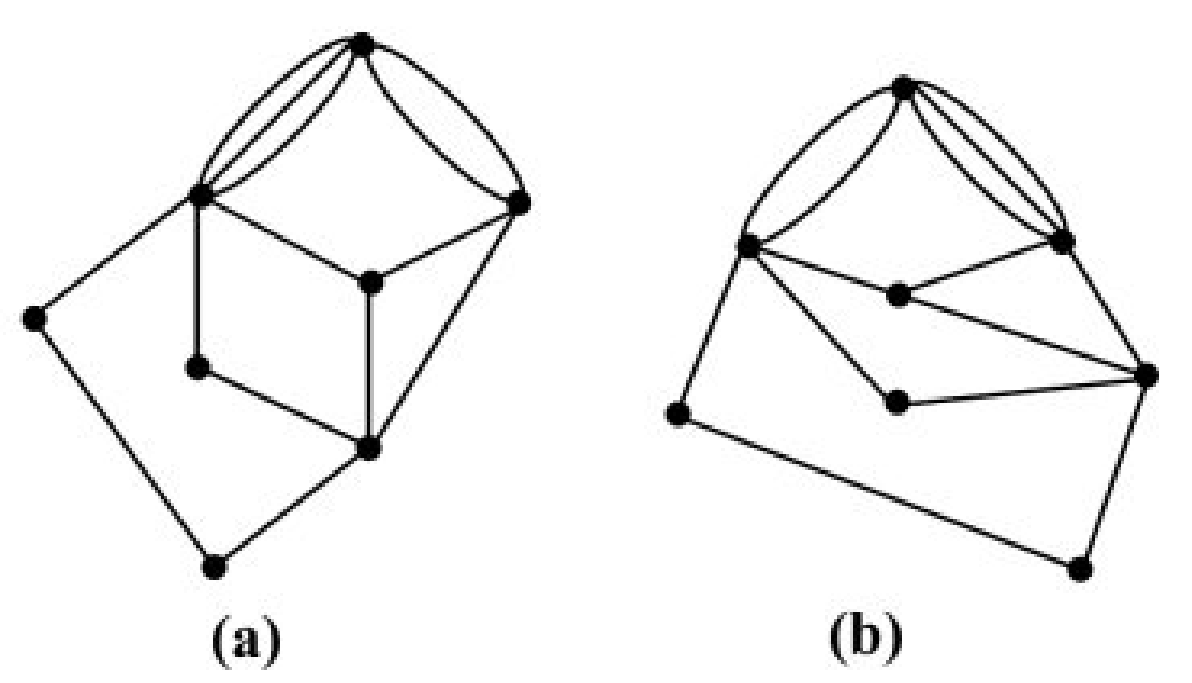,width=2.50in}} \vspace*{8pt}
\caption{(a) Self-dual graph $G(.(2,3).(3,2))$; (b) graph
$G(.(2,3).(2,3))$ which is not self-dual. \label{f1.6}}
\end{figure}

\begin{definition}
An amphicheiral alternating knot or link $L$ is called {\it Dasbach-Hougardy
link} if it has no minimal diagram for which $G(D)$ is isomorphic to
$G^*(D)$.
\end{definition}

According to a computer search, the smallest Dasbach-Hougardy link is 12-crossing alternating link
$(2\,1,2)\,1\,1\,(2\,1,2)$ which belongs to the same family
$(p\,1,2)\,1\,1\,(p\,1,2)$ ($p\ge 2$) as Dasbach-Hougardy counterexample
$(2\,1,3)\,1\,1\,(2\,1,3)$. In this way it is possible to
obtain an infinite number of Dasbach-Hougardy links belonging to the same
family.

Moreover, we propose the more general construction of Dasbach-Hougardy links:

\begin{definition}
Alternating pretzel (Montesinos) tangle $p_1,p_2,\ldots ,p_n$, where
$p_1$, $p_2$, ... $p_n$ ($n\ge 2$) are rational tangles not
beginning by 1 is called {\it oriented} if it is not equal to its
reverse. If all $p_i$ are integers, it is called {\it integer
tangle}, and if at least one $p_i$ is not an integer it is called
{\it non-integer tangle}.
\end{definition}

The symbol $1^{4k-2}$ denotes $1\,1\,...\,1$, where $1$ occurs
$4k-2$ times ($k\ge 1$).
\bigbreak

\noindent {\bf Conjecture 1.2.} {\it Every link given by Conway
symbol of the form
$$(p_1,p_2,\ldots ,p_n)\,1^{4k-2}(p_1,p_2,\ldots ,p_n)$$

\noindent ($k\ge 1$, $n\ge 2$), where $p_1,p_2,\ldots ,p_n$ is
oriented non-integer tangle is Dasbach-Hougardy link. All knots or links of the form

$$(p_1,p_2,\ldots ,p_n)\,t\,(p_1,p_2,\ldots ,p_n)$$

\noindent where $p_1,p_2,\ldots ,p_n$ is oriented integer or
non-integer tangle and $t$ is palindromic\footnote {A rational
tangle is called palindromic if it is equal to its reverse.}
rational tangle are amphicheiral and satisfy the original Kauffman
conjecture.}

The proposed general construction gives an infinite class of Dasbach-Hougardy
links. For example, $(2\,1,2,2)\,1\,1\,(2\,1,2,2)$,
$(3\,1,2,2\,1)\,1\,1\,(3\,1,2,2\,1)$ (Figs. 10, 11), or
$(2\,1,2,2)\,1\,1$ $1\,1\,1\,1\, (2\,1,2,2)$ are examples of Dasbach-Hougardy
links. Each of them can be used as the counterexample to the
original Kauffman conjecture.

\begin{figure}[th]
\centerline{\psfig{file=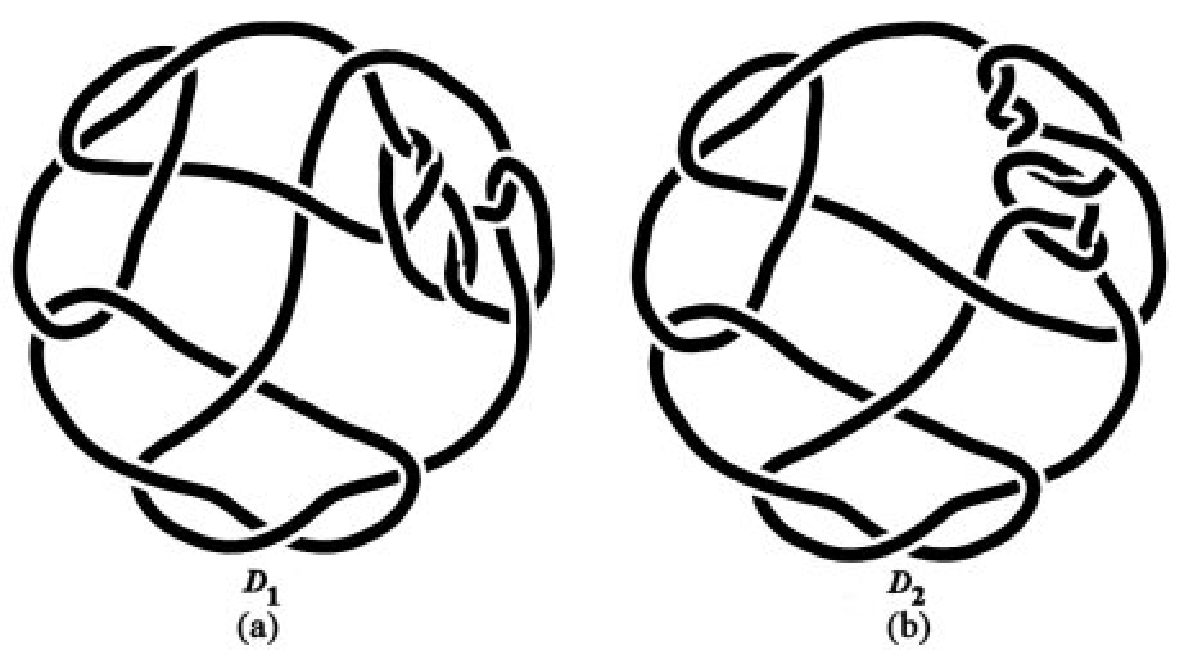,width=2.40in}} \vspace*{8pt}
\caption{Minimal diagrams (a) $D_1$; (b) $D_2$ of the amphicheiral link
$(3\,1,2,2\,1)\,1\,1\,(3\,1,2,2\,1)$. \label{f1.7}}
\end{figure}

\begin{figure}[th]
\centerline{\psfig{file=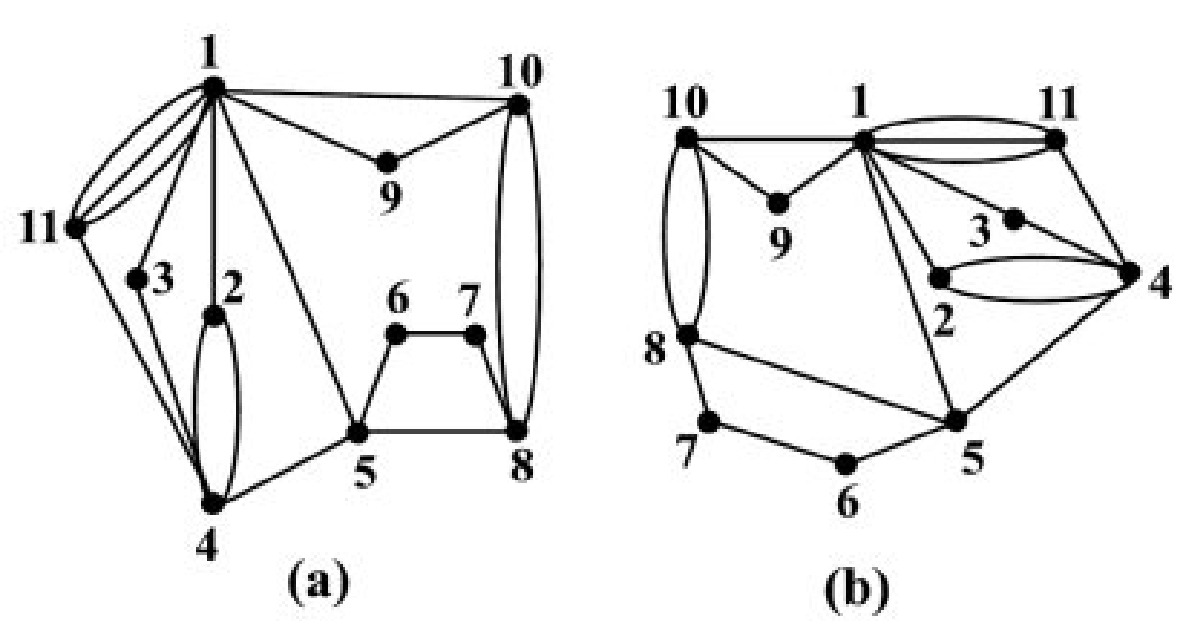,width=3.00in}} \vspace*{8pt}
\caption{Isomorphic graphs $G_1\simeq G^*_2$ of the amphicheiral link $(3\,1,2,2\,1)\,1\,1\,(3\,1,2,2\,1)$. \label{f1.8}}
\end{figure}

All knots or links of the form
$$(p_1,p_2,\ldots ,p_n)\,t\,(p_n,\ldots ,p_2,p_1)$$

\noindent where $p_n,\ldots ,p_2,p_1$ is reverse of $p_1,p_2,\ldots
,p_n$, and $t$ is a palindromic rational tangle are amphicheiral and
satisfy original Kauffman conjecture.

In the paper \cite{3} the authors proved that the original Kauffman Conjecture
is true for negative amphicheiral alternating knots. They also announced the counterexample
to the conjecture that every Dasbach-Hougardy knot is algebraic,
which we proposed in the preceding version of this paper
(arXiv:1005.3612v1 [math.GT]).


\begin{thebibliography}{9}

\bibitem{1} Conway, J. (1970) An enumeration of knots and links and some of
their related properties, in {\it Computational Problems in Abstract
Algebra}, Proc. Conf. Oxford 1967 (Ed. J. Leech), 329--358, Pergamon
Press, New York.

\bibitem{2} Dasbach, O.T., Hougardy, S. (1996) A conjecture of Kauffman on
amphicheiral alternating knots, {\it J. Knot Theory Ramifications},
{\bf 5}, 629--635.


\bibitem{3} Ermotti, N., van Quach Hongler, C., Weber, C.  (2011)
A proof of Tait's Conjecture on alternating $-$achiral knots,
arXiv:1103.3203v1 [math.GT]

\bibitem{4} Jablan, S.V., Sazdanovi\' c, R. (2007) LinKnot- Knot Theory by
Computer, World Scientific, New Jersey, London, Singapore.

\bibitem{5} Kaufman, L.H. (1987) State Models and the Jones Polynomial,
{\it Topology} 26 (1987), 395--407.

\bibitem{6} Kauffman, L.H. (1990) {\it Problems in Knot Theory} in
van Mill, J. and G. Reed, G.M. (Editors) {\it Open Problems in
Topology}, 487--522, Elsevier Science Publishers B.V., North
Holland.

\bibitem{7} van Mill, J., G. Reed, G.M. (1991) {\it Open
Problems in Topology}, Topology Appl. {\bf 42}, 301--307.

\bibitem{8} Menasco, W.W., Thistlethwaite, M.B. (1991)
The Tait flyping conjecture. {\it Bull. Amer. Math. Soc.} (N.S.) 25, {\bf
2}, 403--412.

\bibitem{9} Whitney, H. (1933) 2-Isomorphic Graphs. {\it Amer. J. Math.},
{\bf 55}, 1-4, 245--254.


\end{thebibliography}
\end{document}